\magnification\magstep1
\baselineskip=18pt 
\let\A=\i
\def\w{\widetilde}
\def\i{\infty}
\def\v{\vert}
\def\V{\Vert}

\def\u{{\cal{U}}}
\magnification\magstep1   
\baselineskip=18pt

\def\i{\infty}
\def\l{\ell}
 
\def\conv{{\rm conv}}

\def\ra{\rightarrow}
 
\def\w{\widetilde}
\def\i{\infty}
\def\v{\vert}
\def\V{\Vert}
\def\n{\noindent}
\def\N{{\bf N}}
\def\qed{{\vrule height7pt width7pt
depth0pt}\par\bigskip}
\def\e{\epsilon}
\centerline{\bf  The
$K_t-$functional for the interpolation couple
$L_1(A_0),L_\i(A_1)$} \vskip12pt \centerline{{by 
Gilles Pisier}\footnote*{Supported in part by
N.S.F. grant DMS 9003550}} \vskip12pt

{\bf Abstract}  Let $(A_0,A_1)$ be a compatible
 couple of Banach spaces in the interpolation
theory sense. We give a formula for the
$K_t-$functional of the  interpolation couples
$(\l_1(A_0),c_0(A_1))$ or $(\l_1(A_0),\l_\i(A_1))$ and 
$(L_1(A_0),L_\i(A_1))$.

\vskip12pt
We first recall  the definition of the $K_t$  
functional which is a fundamental tool in the
Lions-Peetre Interpolation Theory   
and also in Approximation Theory, cf.e.g. [1,2].
Let $(A_0,A_1)$ be a compatible couple
of Banach (or quasi-Banach) spaces. This just means that
$A_0,A_1$ are continuously included into a larger
topological vector space (most of the time left implicit),
so that we can consider unambiguously the sets $A_0+A_1$
and $A_0\cap A_1$.  For all $x\in A_0+A_1$ and
for all $t>0$, we let $$K_t(x;A_0,A_1)= \inf
\big({\|x_0\|_{A_0}+t\|x_1\|_{A_1}\ | \ x=x_0+x_1,x_0\in
A_0,x_1\in A_1}).$$  
Recall that the (real interpolation) space
$(A_0,A_1)_{\theta,p}$ is defined ($0<\theta<1,1\leq
p\leq\i$) as the space of all $x$ in $A_0+A_1$ such that
$\|x\|_{\theta,p} <\infty$ where
 $$\|x\|_{\theta,p} =(\int{(t^{-\theta}K_t(x;A_0,A_1))^{p}
dt/t})^{1/p} .$$
It is well known that the $K_t$ functional for the couple
$(L_1(\mu),L_\i(\mu))$ on a non-atomic measure space
$(\Omega,\mu)$ is given by
$$K_t(f;L_1(\mu),L_\i(\mu))=\sup\{\int_E \v f\v d\mu ,\
E\subset \Omega,\  \mu(E)\le t \}.$$
Let $(\w{\Omega},\w{\mu})$ be the measure space obtained
by forming the disjoint union of a sequence of copies of
$(\Omega,\mu)$. Since $L_p(\Omega,\mu;\l_p)$ can be
identified with $L_p(\w{\Omega},\w{\mu})$, we have, for
all $f=(f_i)$ in  $L_1( \mu;\l_1)+L_\i( \mu;\l_\i)$ 
$$\eqalign{ K_t(f;L_1( \mu;\l_1),L_\i( \mu;\l_\i)) 
&=\sup\{ \sum \int_{E_i}\v f_i\v d\mu ,\ \ E_i\subset
\Omega, \ \sum\mu(E_i)\le t\}\cr
&=
\sup\{\sum K_{t_i}(f_i;L_1(\mu),L_\i(\mu)),\ \ t_i\ge 0,\
\sum t_i \le t \}  .}$$

\n Since $L_p( \mu;\l_p)$ and $\l_p(L_p( \mu))$ can be
identified, this example is the prime motivation
for the following statement.

\proclaim Theorem 1. Let $(A_0,A_1)$ be a compatible couple
of Banach spaces. Consider the pair
$(\l_1(A_0),\l_\i(A_1))$. Then,    $\forall x=(x_i)
\in \l_1(A_0)+\l_\i(A_1)$, if $x_i=0$ except for finitely
many indices, we have
$$K_t(x;\l_1(A_0),\l_\i(A_1))=\sup\{\sum_i
K_{t_i}(x_i;A_0,A_1),\ t_i\ge 0,\ \sum t_i\le
t\}.\leqno(1)$$ As a consequence, $\forall x=(x_i) \in
\l_1(A_0)+c_0(A_1)$, we have
$$K_t(x;\l_1(A_0),c_0(A_1))=\sup\{\sum_i
K_{t_i}(x_i;A_0,A_1),\ t_i\ge 0,\ \sum t_i\le t\}. $$

\n{\bf Proof:} Let us denote by $C_t$ the right hand side
of the above identity (1). Then it is very easy to check
that $C_t\le K_t(x;\l_1(A_0),\l_\i(A_1))$. Let us check 
the converse. Let $x$ be
such that $C_t<1$.  This means 
$$\sup_{\sum t_i\leq t}\big\{  \inf_{x_i=a_i+b_i}\big( \sum
\V  a_i\V _{A_0} + t_i\V b_i\V_{A_1} \big)\big\} <1
.\leqno(2)$$ 
We want to deduce from this the same inequality  
but with the inf and the sup interchanged. This can be
viewed as a consequence of the minimax lemma (which
itself is an application of the Hahn-Banach theorem). We
prefer to deduce it directly from the Hahn-Banach
theorem, as follows. This inequality (2) clearly  implies
(choosing $t_i= t\xi_i $) that for any
non-negative sequence $\xi=(\xi_i)$ such that $\sum
\xi_i <1$ there is, for each index $i$  a decomposition
$x_i=\alpha_i+\beta_i$ in $A_0+A_1$ such that $$\sum_i 
\xi_i [(\sum_k\V \alpha_k\V_{A_0}) + t\V
\beta_i\V_{A_1}] <1.\leqno(3)$$ Fix a number $\epsilon
>0$. We will show that the left side of (1) is
  less than $1+\epsilon$. We   assume that, for some
$n$, we have  $x_i=0$ for all indices $i\ge n$. Let
$C\subset {\bf R}^n$ be   the set of all 
points $y=(y_i)$ of the form $$y_i=(\sum_{k\ge 0} \V
a_k\V_{A_0}) + t\V b_i\V_{A_1}\ \hbox{where}\ 
x_i=a_i+b_i,\ a_i\in A_0,\ b_i\in A_1.$$  We claim that
the convex hull of $C $, denoted by $\conv(C)$  intersects
$]-\i,1+\epsilon[^n$. Otherwise,  by Hahn-Banach (we
separate a convex set from an open convex one) we would
find a separating functional $\xi$  and a real number
$r$ such that $\xi <r$ on $]-\i,+1]^n$ and $\xi>r$ on $C$.
But (since we oviously can assume $r=1$) this would
contradict (3). This   shows that $\conv(C)$ intersects
 $]-\i,1+\epsilon[^n$, hence we can find
decompositions $x_i=a_i^m+b_i^m$, $1\le m\le M$ and
positive scalars  $\lambda_1,...,\lambda_m,...,\lambda_M$
with $\sum_m \lambda_m =1$, such that we have for every
index $i$
$$\sum_m \lambda_m [(\sum_{k\ge 0} \V
a_k^m\V_{A_0}) + t\V b_i^m\V_{A_1}]\le
1+\epsilon.\leqno(4)$$ We can then set 
$$ a_i= \sum_m \lambda_m a_i^m,\ b_i=\sum_m \lambda_m
b_i^m.$$
Note that $x_i=a_i+b_i$.
Moreover, by (4) and the triangle inequality,
for every index $i$ $$\sum_{k\ge 0} \V a_k \V_{A_0}+ t
\V b_i\V_{A_1} \le 1+\epsilon ,$$
which clearly implies
$K_t(x;\l_1(A_0),\l_\i(A_1))\le 1+\epsilon$. By
homogeneity, this completes the proof of (1), and
the last assertion is immediate.
 \qed

I asked B.Maurey for some help to extend the preceding
statement without unpleasant asumptions and he
kindly pointed out to me the following fact and its proof: 

\proclaim Theorem 2. Let $P_n$ denote the
projection from $\l_1(A_0)+\l_\i(A_1)$ onto
$\l_1(A_0)+\l_\i(A_1)$ which preserves the first $n$
coordinates and annihilates the other ones.  Then
$$\forall x  \in \l_1(A_0)+\l_\i(A_1)\quad 
 K_t(x;\l_1(A_0),\l_\i(A_1))=\sup_n K_t(P_n
(x);\l_1(A_0),\l_\i(A_1)).\leqno(5)$$ 

\n{\bf Proof:} Fix $t>0$. Clearly the right hand side of
(5) is not more than its left hand side. Conversely, assume
that the right hand side of (5) is $<1$. We will show that
the left side also is less than $1$. To clarify the
notation, if $x$ is a sequence of elements in a Banach
space, we denote by $x(k)$ the $k-$th coordinate of $x$.
Then, for all $x$ as in (5) and for all $n$, there is a
decomposition $P_n(x)=x_0^n + x_1^n$ such that 
$$\V x_0^n\V_{\l_1(A_0)} + t \V x_1^n\V_{\l_\i(A_1)}
<1.\leqno(6)$$  Let $\cal{U}$ be a non trivial
ultrafilter on the positive integers. We   let $n$
tend to infinity along $\u$ and we denote simply by
$\lim_{\u}$ the various resulting limits. Let  $$R=
\lim_{\u} \V x_1^n\V_{\l_\i(A_1)} \quad {\rm
and}\quad a_k=\lim_{\u} \V x_0^n (k)\V_{A_0}.$$ Observe
that (6) implies $$\forall K\in {\N} \quad (\sum_{k<K} a_k
)+tR \leq 1.\leqno(7)$$ Now fix $\e>0$. For each integer
$k$ we can find an integer $n_k>k$ large enough so that 
$$\V x_0^{n_k}(k)\V_{A_0} < a_k +\e 2^{-k} \quad {\rm
and}\quad \V x_1^{n_k} \V_{\l_\i(A_1)} < R+\e.$$
Then we can define
$$  x_0(k) =x_0^{n_k}(k) \quad
{\rm and}\quad x_1(k) =x_1^{n_k}(k) .$$
Clearly $x(k)=x_0(k)+x_1(k)$ for all $k$, and moreover
$$\forall K \sum_{k<K} \V x_0(k)\V_{A_0} +t \sup_{k<K}\V
x_1(k)\V_{A_1} <\sum_{k<K} a_k + \e 2^{-k} +t(R+\e)$$
hence by (7) 
$$\leq 1+ \e (2+t).$$
Since this holds for all $K$, we conclude that $x_0 \in
\l_1(A_0),\ x_1\in \l_\i(A_1)$ and \break $\V
x_0\V_{\l_1(A_0)} +t \V x_1\V_{\l_\i(A_1)} \leq 1+ \e
(2+t),$ and since $\e>0$ is arbitrary we indeed finally
obtain  $$K_t(x;\l_1(A_0),\l_\i(A_1)) \leq 1.{\vrule height7pt width7pt
depth0pt}$$  

\proclaim Corollary 3. The  formula (1)
in the above theorem 1 is valid without any restriction on
$x\in \l_1(A_0)+\l_\i(A_1).$

\n{\bf Remark 4.} The formula (1) remains valid with
the same proof as above if the spaces $A_0$ and $A_1$ are
replaced by families of Banach spaces respectively
$(A_0^n)$ and $(A_1^n)$ . Let us denote by $\l_1(
\{A_0^n\})$ and $\l_\i(\{A_1^n\})$ the corresponding
spaces (these are sometimes called the direct sum of the
families $(A_0^n)$ and $(A_1^n)$ respectively in the
sense of $\l_1$
  and $\l_\i $).This gives us the following generalized
version of (1): 
for all x in $\l_1(
\{A_0^n\})+\l_\i(\{A_1^n\}) $
$$K_t(x;\l_1(\{A_0^n\}),\l_\i(\{A_1^n\}))=\sup\{\sum_i
K_{t_i}(x_i;A_0^i,A_1^i),\ t_i\ge 0,\ \sum t_i\le
t\}.\leqno(8)$$ 

\def\e{\epsilon}
\def\s{{\cal A}}
\def\t{{\cal B}}
 We now reformulate our result in the function space case.
\proclaim Theorem 5. Let $(A_0,A_1)$ be a compatible couple
of Banach spaces. \break Let $(\Omega,\s,\mu)$ be an
arbitrary measure space.  Consider a
function $f$ in \break
$L_1(\Omega,\s,\mu;A_0)+L_\i(\Omega,\s,\mu;A_1)$, where
we define the Banach space valued $L_p$-spaces in the
Bochner sense. Then, for all $t>0$
$$K_t(f;L_1(\Omega,\s,\mu;A_0),L_\i(\Omega,\s,\mu;A_1))=
\sup_{\int \phi d\mu\leq t} \int
K_{\phi(\omega)}(f(\omega);A_0,A_1)
d\mu(\omega),\leqno(9)$$ where the sup runs over all
non-negative measurable functions $\phi$ defined on
$(\Omega,\s)$ with integral not more than $t$.

\n{\bf Proof:} We may clearly assume that the measure
space is $\sigma$-finite.   Now given a function $f_0\in 
L_1(\Omega,\s,\mu;A_0)$, we know
(by definition of Bochner measurability, see e.g. [5]
p.42) that there is a  countable measurable partition of
$\Omega$ into pieces on each of which the oscillation of
$f_0$ for the norm of $A_0$ is small. Similarly, given
$f_1\in L_\i(\Omega,\s,\mu;A_1)$ we know that there is a
measurable partition of $\Omega$ into pieces on each of
which the oscillation of $f_1$ for the norm of $A_1$ is
small.  On the other hand, since the measure space is
$\sigma$-finite, it admits a countable  measurable
partition into sets of finite measure, so that, by
refining the partitions, we can always assume that the
sets have finite measure (so that the conditional
expectation makes sense) and that the same partition
works for both $f_0$ and $f_1$. Consequently, for each
$\e>0$ there is a  countable measurable partition of
$\Omega$ into sets of finite measure on each of which
both  the $A_0$-oscillation of $f_0$  and
 the $A_1$-oscillation of $f_1$
are less than $\e$. The point of this discussion is the
following. Given $f\in
L_1(\Omega,\s,\mu;A_0)+L_\i(\Omega,\s,\mu;A_1)$, we can
find a  $\sigma$-subalgebra $\t\subset \s$ generated by a
countable  measurable partition  of $\Omega$ into sets of
finite measure    such that, if we denote by $f^\t$ the
conditional expectation of $f$ with respect to $\t$, we
have
$$K_t(f-f^\t;L_1(\Omega,\s,\mu;A_0),L_\i(\Omega,\s,\mu;A_1))
<\e.$$ This reduces the proof of (9) to the case when $\s$
is generated by a countable  measurable partition  of
$\Omega$ into sets of finite measure. In that case, we can
identify
 $L_1(\Omega,\s,\mu;A_0)$ and $L_\i(\Omega,\s,\mu;A_1))$
with suitable sequence spaces and (9) follows easily from
(8), (by incorporating the weight of each set of the
partition into the norm of the corresponding
coordinate).\qed

In the situation of Theorem 5, let us assume (for
simplicity) that the intersection $A_0\cap A_1$  is dense
in $A_0$. Then (cf.[1] p.303) we can write for all $x\in
A_0+A_1$
$$ K_t(x;A_0,A_1)=\int_0^{t} k(x,s;A_0,A_1) ds,$$
where the $k-$functional $k(x,s;A_0,A_1)$ is a uniquely
defined nonnegative, nonincreasing, right-continuous
function of $s>0$. In the case of the   (scalar
valued) couple $(L_1,L_\i)$ over a $\sigma$-finite
measure space, we find (cf.[1] p.302)
$$k(x,s;L_1,L_\i)=x^*(s)$$ where $x^*$ is the
nonincreasing rearrangement of $|x|$.

\n Recall the notation
$x^{**}(t)=t^{-1}\int_0^t x^*(s)ds$, so that
$K_t(x;L_1,L_\i)=t x^{**}(t)$. If $0<p\leq
\i,1\leq q\leq\i$ we also recall the 
definition of the  quasi-norm  $\V x\V_{p,q}$
in the Lorentz space
$L_{p,q}$ over a $\sigma$-finite measure space as follows
$$\V x\V_{p,q}=\big
(\int_0^\i[t^{1/p}x^{*}(t)]^q {{dt} \over t} \big
)^{1/q}$$ with the usual convention when $q=\i$.

If $1<p\leq
\i,1\leq q\leq\i$, then Hardy's classical inequality
shows that this is equivalent to the following norm
   $$\V x\V_{(p,q)}=\big
(\int_0^\i[t^{1/p}x^{**}(t)]^q {{dt} \over t} \big
)^{1/q}$$ with the usual convention when $q=\i$. In
particular $L_{p,p}$ is the same as $L_p$ with an
equivalent norm.

With this notation,
we can state

\proclaim Corollary 6. In the same situation as Theorem 5,
assuming (for simplicity) that the intersection
$A_0\cap A_1$  is dense in $A_0$, we denote  for all
$f$ in  
$L_1(\Omega,\s,\mu;A_0)+L_\i(\Omega,\s,\mu;A_1)$, 
$$\forall
s>0\ \forall\omega
\in
\Omega\quad\Psi_f(s,\omega)=k(f(\omega),s;A_0,A_1).$$
Then we have
$$K_t(f;L_1(\Omega,\mu;A_0),L_\i(\Omega,\mu;A_1))=
K_t(\Psi_f;L_1(\Omega\times ]0,\i[, d\mu
ds),L_\i(\Omega\times ]0,\i[, d\mu
ds)).\leqno(10)    $$
Moreover, for $1<p\leq \i,1\leq
q\leq\i$ and $1/p=1-\theta$, we have
$$\V
f\V_{(L_1(\Omega,\mu;A_0),L_\i(\Omega,\mu;A_1))_{\theta,q}}
=\V \Psi_f \V_{(p,q)}\leqno(11)$$
where the     Lorentz space norm is relative to
the product space $(\Omega\times ]0,\i[, d\mu ds)$.

\n{\bf Proof:} By (9) we have
$$K_t(\Psi_f;L_1(\Omega\times ]0,\i[, d\mu
ds),L_\i(\Omega\times ]0,\i[, d\mu
ds))$$
$$=\sup_{\int \phi d\mu\leq t} \int K_{\phi(\omega)}
(\Psi_f(.,\omega);L_1(]0,\i[,ds),L_\i(]0,\i[,ds))
d\mu(\omega)$$
using (9) again this yields (10) since we have obviously 
$$\forall t>0,\forall \omega\in \Omega \quad 
K_t(\Psi_f(.,\omega);L_1(]0,\i[,ds),L_\i(]0,\i[,ds))$$
$$=
\int_0^t \Psi_f(s,\omega)
ds 
 =K_{t}(f(\omega);A_0,A_1).$$  Clearly (11) is an
immediate consequence of (10) by applying 
$K_t(x;L_1,L_\i)=tx^{**}(t)$  on the product space with
$x=\Psi_f$ .\qed

\n{\bf Remark 7. } As an application of Corollary 6, one
can   derive the well known Lions-Peetre results on
interpolation between vector valued $L_p$-spaces in a
rather transparent way, for example  
  in the situation of Corollary 6, if $q=p$ and 
$1/p=1-\theta$, we have
$$(L_1(\Omega,\s,\mu;A_0),L_\i(\Omega,\s,\mu;A_1))_{\theta,p}
=L_p(\Omega,\s,\mu;(A_0,A_1)_{\theta,p}).$$
Indeed, when $p=q>1$ Hardy's classical inequality (see
[1] p.124 and 219) shows that
for all $x$ in $A_0+A_1$, $\V k(x,s;A_0,A_1)\V_{L_p(ds)}$
is equivalent to the norm of $x$ in
$(A_0,A_1)_{\theta,p}$. Therefore, since   $\V \Psi_f
\V_{(p,p)}$ is equivalent to $\V \Psi_f
\V_{L_p(d\mu ds)}$, it is equivalent to the norm  of $f$
in $L_p(\Omega,\s,\mu;(A_0,A_1)_{\theta,p}).$

\n In fact, one finds more generally that if
$1/p=1-\theta$ then for all $1\leq q\leq p$ the following
well known
inclusion holds
$$(L_1(\Omega,\s,\mu;A_0),L_\i(\Omega,\s,\mu;A_1))_{\theta,q}
\subset L_p(\Omega,\s,\mu;(A_0,A_1)_{\theta,q}).$$
Moreover, when $q\geq p$ the reverse inclusion holds. We
refer to [4] for counterexamples to the other inclusions.

\n {\bf Remarks.} 
(i) Using the "power theorem" (cf.[2] p.68) it is easy to
deduce from Theorem 5 an equivalent of the $K_t$-functional
for the couple 
$L_p(\Omega,\s,\mu;A_0),L_\i(\Omega,\s,\mu;A_1)$ for
$0<p<\i$, when $(A_0,A_1)$ are Banach spaces.

\n (ii) More generally, if $1\leq p_0,p_1<\i$ then there
are simple natural quantities known to be equivalent to the
$K_t$-functional for the couple
$L_{p_0}(\Omega,\s,\mu;A_0),L_{p_1}(\Omega,\s,\mu;A_1)$.
In the case $p_1$ finite, these can be derived easily
from the trivial case $p_0=p_1$ and the power theorem,
and this argument even works when $(A_0,A_1)$ are
quasi-Banach spaces. This application of the power theorem
was pointed out to me by Quanhua Xu, but Cwikel informed
me that this was already known to J.Peetre, (cf.also [8]).
Apparently however this approach does not yield the case
$p_1=\i$ which is the main point of the present paper.

\def\o{\omega}
We will give as an application  a generalization of an
embedding theorem for $L_p$ spaces, namely the following.
If $(\Omega',\s',\mu')$  is an arbitrary
measure space, we can define a linear operator
$$T_p:L_p(\Omega' ,\mu')\ra L_{p,\i}(\Omega'\times ]0,\i[
,d\mu' ds) $$ as follows (here $0<p<\i$ and   we
intentionally denote below by $\omega$ a positive real
number instead of $s$ and change the notation $ds$ to
$d\omega$) $$\forall f\in L_p(\Omega' ,\mu') \quad
T_p(f)(\omega',\omega)=\omega^{-1/p}f(\omega').$$ Then it
is a simple exercise to check that $T_p$ is an isometric 
embedding i.e. we have  $$\forall f\in L_p(\Omega' ,\mu')
\quad \V T_p(f)\V_{p,\i}= \V f \V_p.\leqno(12)$$
Actually, if we denote by $m$ the product measure
$dm=d\mu'\times d\omega$, we have 
$$\forall t>0 \quad t^pm(\{|T_p(f)|>t\}) =\int |f|^p
d\mu'.\leqno(13)$$
Similarly, let us denote by $\nu$ the counting
measure on the set $\N^*$ of all positive integers. Then
the preceding embedding has the following discrete
counterpart. We define a linear operator
$$S_p:L_p(\Omega' ,\mu')\ra
L_{p,\i}(\Omega'\times \N ^* ,d\mu' d\nu) $$
as follows ( $0<p<\i$)
$$\forall f\in L_p(\Omega' ,\mu') \quad
S_p(f)(\omega',n)=n^{-1/p}f(\omega').$$
Again, it is easy to check that
$$\forall f\in L_p(\Omega' ,\mu') \quad
\V S_p(f)\V_{p,\i}= \V f \V_p.$$
Moreover, if we denote, for any  positive
real $r$,  
 by $[r]$ the largest integer $n<r$, and if we denote by
$m'$ the product measure $dm'=d\mu'\times d\nu$, we
clearly have $$\forall t>0 \quad  m'(\{|S_p(f)|>t\})
=\int [{{|f|^p}\over{ t^p}}] d\mu'. $$
We now return to the abstract case

\proclaim Theorem 8.  In the same situation as Theorem 5,
assuming (for simplicity) that the intersection
$A_0\cap A_1$  is dense in $A_0$, we define more
generally two linear operators
$$T_p:(A_0,A_1)_{\theta,p}\ra
(L_1(]0,\i[,d\omega;A_0),L_\i(]0,\i[,d\omega;A_1))_{\theta,\i}$$
$$S_p:(A_0,A_1)_{\theta,p}\ra
(L_1(\N ^* ,d\nu;A_0),L_\i(\N ^*,d\nu;A_1))_{\theta,\i}=
(\l_1( A_0),\l_\i( A_1))_{\theta,\i}$$
by setting
$$\forall
x\in (A_0,A_1)_{\theta,p}\quad
  T_p(x)=(\omega\ra \omega^{-1/p}x) \quad {\rm and }\quad
S_p(x)=(n\ra n^{-1/p} x).$$
Then we have $\forall
x\in (A_0,A_1)_{\theta,p}$
$$ \V T_p(x)\V_
{(L_1(]0,\i[,d\o;A_0),L_\i(]0,\i[,d\o;A_1))_{\theta,\i}}=p' 
 {\Big(\int_0^\i {k(x,s;A_0,A_1)}^p
ds\Big)}^{1/p}.\leqno(14)$$ Therefore, (by Hardy's
inequality) $T_p$ is an isomorphic embedding. Similarly,
$S_p$ is an isomorphic embedding.

\n{\bf Proof: } 
Let $f(\omega)=\omega^{-1/p}x$. Then we have
$$\Psi_f(s,\omega)=\omega^{-1/p} k(x,s;A_0,A_1).$$
Note that by (13) we have $$\forall
t>0\quad \Psi_f^*(t)=t^{-1/p} {\Big(\int_0^\i
{k(x,s;A_0,A_1)}^p ds\Big)}^{1/p}.$$
Hence $\Psi_f^{**}(t)=p't^{-1/p} {\Big(\int_0^\i
{k(x,s;A_0,A_1)}^p ds\Big)}^{1/p}$
and (14) follows from (11) with $q=\i$. The discrete case
is now easy and left to the reader.   \qed  

\n {\bf Remark 9.} We do not see how to completely extend
the preceding facts in the case of quasi-Banach spaces
$A_0,A_1$, with $r<1$ and with $L_r(A_0)$ instead of
$L_1(A_0)$. However, the easy direction in theorems 1 or
5 obviously extends up to a constant. For instance, there
is a constant $c$ such that $\forall x\in
\l_r(A_0)+\l_\i(A_1)$ and $\forall t>0$
$$\sup_{\sum t_i^r \leq t^r}(\sum K_{t_i}
(x_i;A_0,A_1)^r)^{1/r} \leq c
K_t(x;\l_r(A_0),\l_\i(A_1)).\leqno(15)$$

To illustrate the possible uses of theorem 8, we conclude
by an application  to the complex interpolation method
which develops in a more abstract way an idea presented
in [9] in the context of $H^p$ spaces.
Again, let $(A_0,A_1)$ be a compatible
 couple of Banach spaces included in a topological vector
space $V$. Assume moreover that there is a quasi-Banach
space $B$ also included in $V$ and such that for some
$ 0<a<1$ we have $$A_0=(B,A_1)_{a,1}.$$ Let $r=1-a$.
As a typical example of this situation the reader should
think of $B=L_r,A_0=L_1,A_1=L_\i.$
For any $x\in A_0+A_1$, we denote by $S^0(x)$ the sequence
$({x\over n})_{n>0}$ and more generally for any complex
number $z$ we denote by $S^z(x)$ the sequence
$({x\over {n^{1-z}}})_{n>0}$.   Moreover we make the rather
restrictive assumption that $S^0$ defines a bounded
operator from $A_0$ into $(\l_r(B),\l_\i(A_1))_{a,\i}$.
The reader will easily check (as in (12) and (13) above)
that this holds for the preceding example with $B=L_r$.
Then we claim that there is a bounded inclusion mapping
$$\forall 0<\theta <1\quad (A_0,A_1)_\theta \subset
(A_0,A_1)_{\theta,p}\quad{\rm if}\quad {1\over
p}=1-\theta.\leqno(16)$$ See [7] for a somewhat related
result. Let us sketch the proof of (16).  Consider an
element $x$ in the open unit ball of the space
$(A_0,A_1)_\theta$. Then there is an analytic function   
$f$  with values in $A_0+A_1$ on the strip $0<\Re (z)<1$,
which is continuous in the closed strip, such that
$f(\theta)=x$ and such that for all real number $t$,
$f(it)$ is in the   unit ball of $A_0$ and
 $f(1+it)$ is in the   unit ball of $A_1$ (and their
respective norms tend to zero when $t$ tends to
infinity). We now apply Stein's interpolation
principle to the analytic family of operators
$S^z$. Consider  $g(z)=S^zf(z)$. Note that
$g(\theta)=S_p(x)$. For simplicity, let us denote
$C=(\l_r(B),\l_\i(A_1))_{a,\i}$. By our restrictive
assumption we have $\sup_t \V g(it)\V_{C
 } \leq c_0$\break  (where $c_0$,$c_1$,$c_2,$ etc... are
constants)  and trivially we have  
$\sup_t \V g(1+it)\V_{\l_\i(A_1)} \leq 1$. Therefore, we
obtain   $\V g(\theta)\V_{(C,\l_\i(A_1))_\theta} \leq c_1$
 . Since  
$(C,\l_\i(A_1))_\theta\subset
(C,\l_\i(A_1))_{\theta,\i}$, we deduce from the
reiteration principle (cf.[2] p.48) that if
$b=(1-\theta)a+\theta$ we have $\V
g(\theta)\V_{(\l_r(B),\l_\i(A_1))_{b,\i}} \leq c_2$. 
By remark 9 and the same computations as above we
have 
$$\V x\V_{(A_0,A_1)_{\theta,p}} \leq c_3\V
S_p(x)\V_{(\l_r(B),\l_\i(A_1))_{b,\i}},$$ so that
(recalling $g(\theta)=S_p(x)$) we finally find 
$\V x\V_{(A_0,A_1)_{\theta,p}} \leq c_4$. This concludes
the proof of the above claim (16). (The reader should
easily fill the  minor technical gaps that we left to avoid
obscuring the idea.) 
\def\R{{\bf R}}
 Now assume given a closed subspace
$S\subset V$ and let
$$S_0 = S\cap A_0,\qquad S_1  = S\cap A_1,\qquad
\beta=S\cap B.$$ \n Let $Q_0 = A_0/S_0$ , $Q_1 = A_1/S_1$
and $Q=B/\beta$ be the associated quotient spaces. Clearly
$(Q_0, Q_1)$ form a compatible couple since there are
natural inclusion maps $$Q_0 \to V/S \quad {\rm and} \quad
Q_1 \to V/S,$$ and similarly $Q \to V/S$.
Obviously, after composition with the quotient mappings in
the above assumption, we get a bounded map from $A_0$ into 
$(Q,Q_1)_{a,1}$, hence (since the latter vanishes on
$S_0$) we have a
bounded map from $Q_0$ into  $(Q,Q_1)_{a,1}$. Similarly, we
find that the same restrictive assumption as
above is satisfied by the quotient spaces and therefore we
conclude that $$\forall 0<\theta <1\quad (Q_0,Q_1)_\theta
\subset (Q_0,Q_1)_{\theta,p}\quad{\rm if}\quad {1\over
p}=1-\theta.\leqno(17)$$

An alternative to the above restrictive assumption is to
assume the following: there is a Banach space $D\subset 
(A_0+A_1)^{\N}$  and a constant $c$ such that
$$\forall x\in A_0,\ \forall t\in \R\quad \V S^{it}x\V_{D}
\leq c  \V x\V_{A_0} ,\leqno (18)$$
and $$(D,\l_\i (A_1))_{\theta,\i} \subset (\l_1(A_0),\l_\i
(A_1))_{\theta,\i}.\leqno(19)$$
Then (16) holds.
Indeed with the same notation as above,  if $\V
x\V_{(A_0,A_1)_\theta} <1$, this gives $\V
g(\theta)\V_{(D,\l_\i(A_1))_\theta} \leq c_1$, hence a
fortiori $\V g(\theta)\V_{(D,\l_\i(A_1))_{\theta,\i}} \leq
c_2$, therefore by (19) \break $\V
S_p(x)\V_{(\l_1(A_0),\l_\i(A_1))_{\theta,\i}} \leq c_3$,
and by theorem 8, finally  $\V x\V_{(A_0,A_1)_{\theta,p}}
\leq c_4$. Theses assumptions (18) and  (19) are slightly
more general than the preceding one but seem less
easy to verify in practise.

In [9], the preceding argument is applied
  in the particular case $A_0=L_1,\ Q_0=L_1/H^1$,
$A_1=L_\i ,\  Q_1=L_\i/H^\i$ to give a new proof that (17)
holds in this case, which is originally due to Peter Jones
[6]. We refer  the reader to [9] for more information on
this topic. Concerning for instance $H^p$-spaces with
several complex variables or Sobolev spaces on $\R^n$ (cf.
Bourgain's recent paper [3]) the preceding remarks show
that whenever the appropriate real interpolation results
hold, the corresponding complex interpolation results will
also hold. Unfortunately, the real interpolation results do
not seem complete enough at the moment to yield the
assumptions needed in the above remarks.

\vfill\eject

\centerline {\bf References}\vskip6pt

\item {1.} C.Bennett and R.Sharpley, Interpolation of
operators.Academic Press,1988.

\item {2.} J.Bergh and J.L\"ofstr\"om, Interpolation
spaces, An introduction, Springer Verlag 1976.

\item {3.} J.Bourgain, Some consequences of Pisier's
approach to interpolation. Preprint.

\item {4.} M.Cwikel, On
$(L^{p_0}(A_0),L^{p_1}(A_1))_{\theta,q}$.
Proc.Amer.Math.Soc.44 (1974) 286-292.

\item {5.} J.Diestel and J.J.Uhl Jr., Vector
measures.\quad\quad\qquad\quad\quad\quad\quad\quad
\qquad\qquad\qquad\qquad \break 
  Math. Surveys 15.  Amer. Math.Soc. Providence, 1977.

\item {6.} P.Jones, $L^\infty$ estimates for the
$\bar{\partial}$-problem in a half plane. Acta Math. 150
(1983) \nobreak{137-152}.

\item {7.} J.Peetre, Sur l'utilisation des suites
inconditionnellement sommables dans la th\'eorie des
espaces d'interpolation, Rend. Sem. Mat. Univ. Padova 46
(1971) 173-190.

\item {8.} L.Persson, Description of some interpolation
spaces in the off diagonal cases,\break in Interpolation
spaces and related topics in Analysis,Proceedings Lund 1983,
\qquad\quad\  \break 
 Springer Lecture Notes in Math.
1070 (1984) 213-230.

\item {9.} G.Pisier, Interpolation between $H^p$ spaces
and non-commutative generalizations I.\break  To appear.

\vskip12pt

Texas A. and M. University

College Station, TX 77843, U. S. A.

and

Universit\'e Paris 6

Equipe d'Analyse, Bo\^\A te 186,
 
75230 Paris Cedex 05, France

\end